\documentclass[12pt]{article}
\usepackage{amssymb,amsmath,amsthm, amsfonts}
\usepackage{graphicx}
\usepackage{epsfig}
\usepackage{tikz}
\def\dref#1{(\ref{#1})}

\theoremstyle{plain}
\newtheorem{theorem}{Theorem}[section]
\newtheorem{lemma}{Lemma}[section]

\theoremstyle{definition}

\numberwithin{equation}{section}

\begin{document}

\title{\large\bf Asymptotic behavior of solutions to a tumor angiogenesis model with chemotaxis--haptotaxis}

\author
{\rm Peter Y.~H.~Pang\\
\it\small Department of Mathematics, National University of Singapore\\
\it\small 10 Lower Kent Ridge Road, Republic of Singapore 119076
\\[3mm]
\rm Yifu Wang\thanks{
Corresponding author. Email: {\tt wangyifu@bit.edu.cn}}\\
\it\small School of Mathematics and Statistics, Beijing
\it\small Institute of Technology\\
\it\small Beijing 100081, People's Republic of China}

\date{}

\maketitle
\noindent
\begin{abstract}
 This  paper studies the following system of differential equations modeling tumor angiogenesis in a bounded smooth domain $\Omega \subset \mathbb{R}^N$ ($N=1,2$):
$$\label{0}
 \left\{\begin{array}{ll}
  p_t=\Delta p-\nabla\cdotp p(\displaystyle\frac \alpha {1+c}\nabla c+\rho\nabla w)+\lambda p(1-p),\,&
x\in \Omega, t>0,\\
 c_t=\Delta c-c-\mu pc,\, &x\in \Omega, t>0,\\
w_t= \gamma p(1-w),\,& x\in \Omega, t>0,
  \end{array}\right.
$$
where $\alpha, \rho, \lambda, \mu$ and $\gamma$  are positive parameters. For any reasonably regular
initial data $(p_0, c_0, w_0)$, we prove the global boundedness ($L^\infty$-norm) of $p$  via  an
iterative method. Furthermore, we investigate the long-time behavior of solutions to the above system under an additional mild condition,
 and improve previously known results. In particular, in the one-dimensional case, we show that the solution $(p,c,w)$ converges to
 $(1,0,1)$ with an explicit exponential rate as time tends to infinity.



\end{abstract}
\noindent {\bf\em Key words:}~Angiogenesis, chemotaxis, haptotaxis,  boundedness, asymptotic behavior

\noindent {\bf\em 2010 Mathematics Subject Classification}:~35A01, 35B35, 35K57, 35Q92,  92C17

\section{Introduction}

As a  physiological process, angiogenesis involves the formation of new capillary networks sprouting from a pre-existing
vascular network and  plays  an important role in embryo development, wound healing and tumor growth.
For example, it has been recognized that
 capillary growth through angiogenesis
leads to vascularization of a tumor, providing it with its own dedicated blood supply and
consequently allowing for rapid growth and metastasis.
The process of tumor angiogenesis can be divided into
three main stages (which may be overlapping): (i) changes within existing blood vessels;
(ii) formation of new vessels; and (iii) maturation of new vessels. Over the past decade, a lot of work has been done on
the mathematical modeling of tumor growth; see, for example, \cite{Anderson,Bellomo,Chaplain1,Chaplain2,Mu,Surulescu1,Surulescu2}
and the references cited therein. In particular, the role of angiogenesis in tumor growth has also attracted a great deal of attention;
see, for example, \cite{Anderson1,Chaplain3,Levine,Paweletz,Sleeman} and the
references cited therein. 
For example, in Levine et al.~\cite{Levine}, a system of PDEs
using reinforced random walks was deployed to model the first stage of angiogenesis, in which
chemotactic substances from the tumor combine with the
receptors on the endothelial cell wall to release proteolytic enzymes 
that can degrade the basal membrane of the blood vessels eventually.

In this paper we consider a variation of the model proposed in \cite{Anderson}, namely, 
 \begin{equation}\label{1.1}
 \left\{\begin{array}{ll}
  p_t=\Delta p-\nabla\cdotp p(\displaystyle\frac \alpha {1+c}\nabla c+\rho\nabla w)+\lambda p(1-p),\,&
x\in \Omega, t>0,\\
 c_t=\Delta c-c-\mu pc,\, &x\in \Omega, t>0,\\
w_t= \gamma p(1-w),\,& x\in \Omega, t>0, \\[.15cm]
\displaystyle{\frac {\partial p}{\partial \nu }- p(\displaystyle\frac \alpha {1+c}\frac {\partial c}{\partial \nu }+
\rho\frac {\partial w}{\partial \nu }
)=\frac {\partial c}{\partial \nu }=0},\quad
&x\in \partial\Omega, t>0,\\[.15cm]
p(x,0)=p_0(x), \ c(x,0)=c_0(x),  \ w(x,0)= w_0(x), \quad&
x\in \Omega,\\
  \end{array}\right.
\end{equation}
in a bounded smooth domain $\Omega \subset \mathbb{R}^N (N=1,2)$,  where, in addition to random motion, the existing blood vessels' endothelial cells $p$
migrate in response to the concentration gradient 
of a chemical signal $c$ (called Tumor Angiogenic Factor, or TAF) secreted by tumor cells  as well as
the concentration gradient of non-diffusible glycoprotein fibronectin $w$ produced by the endothelial cells \cite{Morales-Rodrigo}.
The former directed migration is a chemotatic process, whereas the latter is a haptotatic process. In this model, it is assumed that
 the endothelial cells proliferate according to a logistic law, that the spatio-temporal evolution of TAF occurs through 
 diffusion, 
 natural decay and degradation upon binding to the endothelial cells, and that the fibronectin is produced by the endothelial cells and
 degrades upon binding to the endothelial cells.

For the remainder of this paper, we shall assume that the initial data satisfy the following:
\begin{equation}\label{1.4}
 \left\{
 \begin{array}{l}
  (p_0,c_0,w_0)\in (C^{2+\beta}(\overline{\Omega}))^3~\hbox{is nonnegative}~\hbox{for some}~\beta\in (0,1) ~  \hbox{with}~p_0\not\equiv  0,\\[.2cm]
 \displaystyle{\frac {\partial p_0}{\partial \nu }- p_0(\displaystyle\frac \alpha {1+c_0}\frac {\partial c_0}{\partial \nu }+
\rho\frac {\partial w_0}{\partial \nu })=\frac {\partial c_0}{\partial \nu }=0}.
  \end{array}
  \right.
\end{equation}

 The present paper focuses on the global existence and asymptotic behavior of classical solutions to \eqref{1.1}.
Let us look at two subsystems contained in \eqref{1.1}. The first is a Keller--Segel-type chemotaxis system with signal absorption: 
\begin{equation}\label{1.2}
\left\{
\begin{array}{ll}
p_t=\Delta p-\nabla\cdot( p\nabla c)+\lambda p(1-p),\quad &x\in\Omega, t>0,
\\
c_t=\Delta c- pc,\quad &x\in\Omega, t>0.
\\
\end{array}\right.
\end{equation}
It is known that, unlike the standard Keller--Segel model, \eqref{1.2} with $\lambda=0$ possesses
global, bounded classical solutions in two-dimensional bounded convex domains for arbitrarily large initial data;
while in three spatial dimensions, it admits at least global weak solutions which
eventually become smooth and bounded after some waiting time \cite{Tao1}. In the high-dimensional setting, it has been proved that
global bounded classical solutions exist for suitably large $\lambda>0$, while only certain
weak solutions are known to exist  
for arbitrary $\lambda>0$ \cite{Lankeit}.

 Another delicate subsystem of \eqref{1.1} is the haptotaxis-only system 
 obtained by letting $\alpha=0$ in \eqref{1.1}:
 \begin{equation}\label{1.3}
 \left\{\begin{array}{ll}
  p_t=\Delta p-\rho\nabla\cdotp(p\nabla w)+\lambda p(1-p),\,&
x\in \Omega, t>0,\\
w_t= \gamma p(1-w),\,& x\in \Omega, t>0. \\
  \end{array}\right.
\end{equation}
Here, since the quantity $w$ satisfies an ODE without any diffusion, the smoothing effect on the spatial regularity of $w$
during evolution cannot be expected.  
To the best of our knowledge, unlike the study of chemotaxis systems, 
the mathematical literature
on  haptotaxis systems is comparatively thin. Indeed, the literature provides only some results on 
 global solvability in various 
special models, and the detailed description of qualitative properties such as long-time behaviors of solutions
is available only in very particular cases (see, for example, \cite{Corrias,Litcanu,Litcanu2,Marciniak,Tao,TWCPAA,Walker,Winkler1}).

More recently, some results on
global existence and asymptotic behavior for certain chemotaxis--haptotaxis models of cancer invasion
have been obtained (see, for example, \cite{YanLi,PW1,PW2,SSW,Tao2,Taosiam,Wangke}). Particularly, Hillen et al.~\cite{Hillen} have shown the convergence of a cancer invasion model in one-dimensional domains and the result has been subsequently extended to higher dimensions \cite{YanLi,Taosiam,Wangke}.


In \cite{Morales-Rodrigo}, in two spatial dimensions, the authors showed the global existence and
long-time behavior of classical solutions to \eqref{1.1} when the initial data $(p_0, c_0, w_0)$
satisfy 
 either $w_0 > 1$ or $ \|w_0-1\|_{L^\infty(\Omega)}<\delta$
for some  $\delta> 0$ (see Lemma 5.8 of \cite{Morales-Rodrigo}).
Generalizing this result, our first main result establishes that, for any choice of reasonably regular
initial data $(p_0, c_0, w_0)$,  
the $L^\infty$-norm of $p$ is globally bounded. This is done via  an
iterative method. 


\begin{theorem}\label{theorem1.1} Let $\alpha, \rho, \lambda, \mu$ and $\gamma$  be  positive parameters. Then for any
 initial data $(p_0, c_0, w_0)$ satisfying \eqref{1.4}, the problem \eqref{1.1} possesses a unique classical solution $(p,c,w)$ comprising nonnegative functions in $C(\bar{\Omega}\times [0,\infty))\cap C^{2,1}(\bar{\Omega}\times (0,\infty)$ such that  $\|p(\cdot,t)\|_{L^\infty(\Omega)}\leq C$ for all $t>0$.
\end{theorem}

Next, we investigate the asymptotic behavior of solutions  to \eqref{1.1}. Under an additional mild condition
on the initial data $w_0$, we will show that the solution $(p,c,w)$ converges to the spatially homogeneous equilibrium $(1,0,1)$ as time tends to infinity.

\begin{theorem}\label{theorem} Let $\alpha, \rho, \lambda, \mu$ and $\gamma$  be positive parameters, and  suppose that \eqref{1.4} is satisfied and $w_0>1-\displaystyle\frac 1 \rho$.
Then the solution $(p,c,w)\in C(\bar{\Omega}\times [0,\infty))\cap C^{2,1}(\bar{\Omega}\times (0,\infty)$ of  \eqref{1.1} satisfies
\begin{equation}\label{1.5}
\displaystyle\lim_{t\rightarrow \infty}\|p(\cdot,t)-1\|_{L^r(\Omega)}+\|c(\cdot,t)\|_{W^{1,2}(\Omega)}+\|w(\cdot,t)-1\|_{L^r(\Omega)}=0
\end{equation}
for any $r\geq 2$. In particular, if $N=1$, then 
for any $\epsilon\in(0,\min\{\lambda_1,1,\gamma,\lambda\})$ there exists $C(\epsilon)>0$ such that
\begin{equation}\label{1.6}
\displaystyle \|p(\cdot,t)-1\|_{L^\infty(\Omega)} \leq C(\epsilon)e^{-(\min\{\lambda_1,1,\gamma,\lambda\}-\epsilon) t},
\end{equation}
\begin{equation}\label{1.7}
\|c(\cdot,t)\|_{W^{1,2}(\Omega)}\leq C(\epsilon)e^{-(1-\epsilon)t},
\end{equation}
\begin{equation}\label{1.8}
\displaystyle \|w(\cdot,t)-1\|_{W^{1,2}(\Omega)} \leq C(\epsilon)e^{-(\gamma-\epsilon) t},
\end{equation}
where $\lambda_1>0$ is the
first nonzero eigenvalue of $-\Delta$ in $\Omega$ with the homogeneous Neumann boundary
condition.
\end{theorem}


The main mathematical challenge of the full chemotaxis--haptotaxis system is the strong coupling between the migratory cells $p$ and the haptotactic agent $w$.
This strong coupling has an important effect on the spatial regularity of $p$ and $w$.
In fact, the lack of regularization effect in the spatial variable in the $w$-equation
and the presence of $p$ therein demand tedious estimates on the solution. The key ideas behind our results are as follows:

As pointed out in \cite{Taosiam}, the variable transformation  $z:= p e^{-\rho w}$ plays an important
role in the examination of global solvability for the full chemotaxis--haptotaxis model in the two- and higher-dimensional setting.
However, due to the presence of the additional chemotaxis term in our model, this approach is not directly applicable to our problem.
Instead, in the derivation of Theorem 1.1, we introduce the variable transformation $q := p(c + 1)^{-\alpha} e^{-\rho w}$ as in \cite {Morales-Rodrigo}, and thereby ensure that
$q(\cdot,t)$  is bounded in $L^ n(\Omega)$ for any finite $n$ (see Lemma \ref{lemma23}).
It is essential to our approach to derive a bound for $\int_{\Omega}q^{2^{m+1}}+\int^{t+\tau}_t\int_{\Omega} |\nabla q^{2^{m}} |^2$ from the bound of $\int^{t+\tau}_t\int_{\Omega}q^{2^{m}}$ ($m=1,2,\ldots$) by making appropriate use of  \eqref{2.3a}--\eqref{2.4a} in Lemma \ref{lemma22} (see \eqref{2.3} below).

The crucial idea of the proof of Theorem 1.2 is to show
$$
\displaystyle\frac d {dt} F(p(t), w(t))+\frac 12 \displaystyle \int_\Omega \displaystyle\frac {|\nabla p|^2}p+
\frac 12 \displaystyle \int_\Omega  p|\nabla w|^2\leq
C\int_\Omega p|w-1|+C\int_{\Omega} p|\nabla c|^2
 $$
 with
 $$
 F(p,w)=\kappa \int_\Omega |\nabla w|^2+ \int_\Omega p(\ln p-1) +\int_\Omega p(w-1)- \gamma \kappa  \int_\Omega p(w-1)^2
$$
for some $C>0$ (see (3.10)), on the basis of the global boundedness of the $L^\infty$-norm  of $p$  provided by Theorem 1.1 and under the mild assumption that $w_0>1-\frac 1 \rho$.
Furthermore, in the one-dimensional case, with the help of  higher regularity estimates of $z$ with $z=p e^{-\rho w} $ (see (3.29)),  we show that $p(\cdot,t)$ converges to 1 uniformly in
$\Omega$ as $t\rightarrow \infty$. From this, we derive the exponential convergence of solutions as desired.


\section{Proof of Theorem 1.1}

In this section, we first recall  
the following estimates for the heat semigroup $(e^{\tau \Delta})_{\tau\geq 0}$ in $\Omega\subset\mathbb{R}^N$ under
the Neumann boundary condition. We shall omit the proof thereof and refer the interested readers to \cite{Cao,Horstmann,Winkler792}.
\begin{lemma}\label{lemma316} Let $(e^{\tau \Delta})_{\tau\geq 0}$ be the Neumann heat semigroup in $\Omega$ and $\lambda_1>0$ the first nonzero eigenvalue of $-\Delta$ with homogeneous Neumann boundary condition. Then there exist positive constants $k_1,k_2, k_3$ such that:

i) If $1\leq q\leq p\leq \infty$, then for all  $\varphi\in L^{q}(\Omega)$ with $\int_\Omega \varphi =0 $,
$$
\| e^{\tau \Delta}\varphi\|_{L^p(\Omega)}
\leq k_1(1+\tau^{-\frac N2(\frac 1 q-\frac 1p)})e^{-\lambda_1 \tau}\|\varphi\|_{L^q(\Omega)};
$$

ii) If $1\leq q\leq p\leq\infty$, then for all $\varphi\in L^q(\Omega)$,
$$\|\nabla e^{\tau\Delta}\varphi\|_{L^p(\Omega)}\leq k_2(1+\tau^{-\frac12-\frac N2(\frac1q-\frac1p)})e^{-\lambda_1\tau}\|\varphi\|_{L^q(\Omega)};
$$

iii) If $1\leq q\leq p\leq \infty$, then for 
all $\varphi\in C^1(\bar{\Omega};\mathbb{R}^N)$ with $\varphi \cdot \nu=0$ on $\partial\Omega$,
\begin{equation*}\label{2.9}
\| e^{\tau \Delta}\nabla\cdot\varphi\|_{L^p(\Omega)}
\leq
k_3(1+\tau ^{-\frac {1}{2}-\frac N2(\frac 1 q-\frac 1p)}) e^{-\lambda_1 \tau}\|\varphi\|_{L^q(\Omega)}.
\end{equation*}

\end{lemma}

Next, we recall the following  result on local existence and uniqueness of classical solutions to \eqref{1.1} as
well as a convenient extensibility criterion, which follows from
Theorem 3.1, Lemma 5.9 and Theorem 5.1  of \cite{Morales-Rodrigo}.
\begin{lemma}\label{lemma21}\rm{(\cite {Morales-Rodrigo})}~\it
Let $\Omega\subset \mathbb{R}^N$ be a smooth bounded domain.
There exists $T_{max}\in (0, \infty]$ such that the problem \eqref{1.1}
  possesses  a unique classical solution satisfying
  $(p,c,w)\in (C(\bar{\Omega}\times [0,T_{\max}))\cap C^{2,1}(\bar{\Omega}\times (0,T_{\max}))^3$. Moreover, for any $s>N+2$,
\begin{equation}\label{2.1}
\limsup_{t\nearrow T_{\max}}\|p(\cdot,t)\|_{W^{1,s}(\Omega)}\rightarrow\infty
\end{equation}
if $T_{\max}<+\infty$.
\end{lemma}
 From now on,
 let $(p,c,w)$ be the local classical solution of (1.1) on $(0,T_{\max})$ provided by
 Lemma 2.1, 
 and $\tau:=\min\{1,\frac {T_{\max}}6\}$.

 The following basic  but  important properties of the solution to \eqref{1.1} can be directly obtained via
standard arguments.
 \begin{lemma}\label{lemma22}\rm{(\cite{Morales-Rodrigo})}~\it  There exists a positive constant $C$ independent of time such that
\begin{equation}
\begin{array}{l}
 \displaystyle\int_{\Omega} p(\cdot,t)\leq c_1:=\max\{\int_\Omega p_0,|\Omega|\},~
\displaystyle \int^{t}_{0}e^{-2s} \int_{\Omega}p^2ds \leq C~~\hbox {for  all}~~t\in (0, T_{\max}),
\end{array}
\end{equation}
\begin{equation}\label{2.3a}
\displaystyle \int^{t+\tau}_{t}\int_{\Omega} p^2\leq c_1(1+\frac 1 \lambda)~~~\hbox{for all}~~t\in (0, T_{\max}-\tau),
\end{equation}
\begin{equation}\label{2.4a}
c(t)\leq \|c_0\|_{L^\infty(\Omega)}e^{-t},~ \displaystyle \int_\Omega |\nabla c(t)|^2+ \int^{t}_{0}\int_{\Omega}(|\nabla c|^2+ |\Delta c|^2)\leq C
~~\hbox {for  all}~~t\in (0, T_{\max}),
\end{equation}
\begin{equation}\label{2.5a}
0\leq w(t)\leq \max\{\|w_0\|_{L^\infty(\Omega)},1\}~~\hbox {for  all}~~t\in (0, T_{\max}).\end{equation}
\end{lemma}

As the proof of Theorem 1.1 in the one-dimensional case is similar to that for two dimensions, henceforth in this section, we shall focus on the case $N=2$.

First, we shall show that  $p$  remains bounded in $L^ n(\Omega)$ for any finite $n $.
We note that the $L^ n(\Omega)$-bound  in Lemma 3.10 of \cite{Morales-Rodrigo} depends on the time variable.

\begin{lemma}\label{lemma23}  
For any $n\in (1, \infty)$, there exists a positive constant $C(n,\tau)$ independent of $t$, such that
$\|p(\cdot,t)\|_{L^n(\Omega)}\leq C(n,\tau)$ for all $t\in (0, T_{\max})$.
\end{lemma}

\it{Proof.}\rm\quad  
Let $q := p(c + 1)^{-\alpha} e^{-\rho w}$. As in the proof of Lemma 3.10 in \cite{Morales-Rodrigo}, we infer that for any $m=1,2, \ldots$ there exist constants $c(m)>0$ depending upon $m$ and $c_1>0$ such that
\begin{equation}\label{2.2}
\displaystyle \frac{d}{dt}\int_{\Omega}q^{2^m}(c+1)^\alpha e^{\rho w}+\int_{\Omega} |\nabla q^{2^{m-1}} |^2
\leq c(m)(\displaystyle\int_{\Omega} |\Delta c |^2+1)\int_{\Omega}q^{2^m}+ c(m)(\int_{\Omega}q^{2^m})^2+c_1.
\end{equation}
Next, we use induction to  show
\begin{equation}\label{2.3}
\int_{\Omega}q^{2^m}+\int^{t+\tau}_t\int_{\Omega} |\nabla q^{2^{m-1}} |^2\leq C(m).
\end{equation}
Taking $m = 1$ in \eqref{2.2}, we get
\begin{equation}\label{2.4}
\displaystyle \frac{d}{dt}\int_{\Omega}q^{2}(c+1)^\alpha e^{\rho w}+\int_{\Omega} |\nabla q |^2
\leq c(1)(\displaystyle\int_{\Omega} |\Delta c |^2+1 +\int_{\Omega}q^{2})\int_{\Omega}q^2(c+1)^\alpha e^{\rho w}+c_1,
\end{equation}
which implies that for the functions $y(t)=\int_{\Omega}q^{2}(c+1)^\alpha e^{\rho w}$ and $a(t)= c(1)(\int_{\Omega} |\Delta c |^2+1 +\int_{\Omega}q^{2})$, we have
$$\frac {d y}{dt}\leq a(t)y +c_1.$$
On the other hand, for any given $t>\tau$, it follows from (2.3) that there exists some $t_0\in [t-\tau,t]$ such that
$y(t_0)\leq \frac {c_1}{\tau}(1+\frac 1 \lambda)$.
Hence by ODE comparison argument we get
\begin{equation}\label{2.5}
y(t)\leq y(t_0)e^{\int^t_{t_0}a(s)ds} +c_1\int^t_{t_0}e^{\int^t_{s}a(\tau)d\tau} ds\leq c_2.
\end{equation}
In this inequality,  we have taken $t_0=0$ if $t\leq \tau$ and  noticed that
$\int ^{t}_{t-\tau} a(s)ds\leq c_3 $ for all $t<T_{\max}$ by Lemma 2.2. Combining \eqref{2.4} with  \eqref {2.5}, one can see that \eqref{2.3} is indeed valid for $m=1$.

Now, suppose that \eqref{2.3} is valid for an integer $ m+1=k\geq 2$, i.e.,
\begin{equation}\label{2.6}
\int_{\Omega}q^{2^{k-1}}+\int^{t+\tau}_t\int_{\Omega} |\nabla q^{2^{k-2}} |^2\leq C(k).
\end{equation}
By the Gagliardo--Nirenberg inequality  in two dimensions
$$\|z\|^4_{L^4(\Omega)}\leq c_3 \| \nabla z\|^2_{L^2(\Omega)}\|z\|^2_{L^2(\Omega)}+c_4 \|z\|^4_{L^2(\Omega)},
$$
and hence
\begin{equation}\label{2.7}
\int_{\Omega}q^{2^{k}}\leq c_3 \int_{\Omega} |\nabla q^{2^{k-2}} |^2 \int_{\Omega}q^{2^{k-1}}+c_4 (\int_{\Omega}q^{2^{k-1}} )^2.
\end{equation}

Integrating \eqref{2.7} between $t$ and $t+\tau$ and taking \eqref{2.6} into account, we have
 \begin{equation}
\int^{t+\tau}_t\int_{\Omega}q^{2^{k}}\leq c_5(k),
\end{equation}
which  implies that for any $t\geq \tau$, there exists some $t_0\in [t-\tau,t]$ such that
$\int_{\Omega}q^{2^{k}} (t_0)\leq c_6$.
At this point, let $y(t):=\int_{\Omega}q^{2^k}(c+1)^\alpha e^{\rho w}$ and $b(t)= c(k)(\int_{\Omega} |\Delta c |^2+1 +\int_{\Omega}q^{2^k})$. Then \eqref{2.2}  can be rewritten as
$$\frac {d y}{dt}+\int_{\Omega} |\nabla q^{2^{k-1}} |^2\leq b(t)y +c_1.$$
By the argument above, one can obtain
\begin{equation}\label{2.9a}
\int_{\Omega}q^{2^{k}}+\int^{t+\tau}_t\int_{\Omega} |\nabla q^{2^{k-1}} |^2\leq C,
\end{equation}
 and thereby conclude that \eqref{2.3} is  valid for all integers $ m\geq 1$. The proof of  Lemma 2.3 is now complete in view of the boundedness of the weight  $(c+1)^\alpha e^{\rho w}$.

\vskip3mm
To establish a priori estimates of $\|p(\cdot,t)\|_{L^\infty(\Omega)}$, we need some fundamental estimates for the solution of the following problem: 
\begin{equation}
 \left\{\begin{array}{l}
  c_t=\Delta c -c+f,\quad
x\in \Omega, t>0,\\
\displaystyle \frac{\partial c}{\partial \nu}=0,\quad
x\in \partial\Omega, t>0,\\
c(x,0)=c_0(x),\quad x\in \Omega.\\
 \end{array}\right.\label{210}
\end{equation}

\begin{lemma}\rm{(\cite[Lemma 2.2]{LWang})}\label{lemma24}~\it
 Let  $T>0$, $r\in (1,\infty)$. 
Then for each $c_0\in W^{2,r}(\Omega)$  with $\displaystyle\frac{\partial c_0}{\partial\nu}=0$ on $\partial \Omega$ and $f\in L^r(0,T;L^r(\Omega))$, \eqref{210} has a unique solution
$c\in  W^{1,r}(0,T; L^r(\Omega))\cap  L^r(0,T; W^{2,r}(\Omega))$ given by
$$
c(t)=e^{-t}e^{t\Delta}c_0+\int^t_0 e^{-(t-s)}e^{(t-s)\Delta}f(s)ds, \quad t\in [0,T],
$$
where $e^{t\Delta}$ is the semigroup generated by the Neumann Laplacian, and there is $C_r>0$ such that
 \begin{equation}
\int ^t_0\int _\Omega e^{rs}|\Delta c(x,s)|^r dxds\leq C_r \int ^t_0\int _\Omega e^{rs} |f(x,s)|^rdxds + C_r\|v_0\|_{W^{2,r}(\Omega)}.
\end{equation}
\end{lemma}

Now applying these estimates to control the cross-diffusive flux appropriately, we can derive the boundedness of $p$ in $\Omega\times (0, T_{max})$.

\begin{lemma}\label{lemma25}  
There exists a constant $C>0$
independent of $t$ such
that $\|p(\cdot,t)\|_{L^\infty(\Omega)}\leq C$ for all $t\in (0, T_{max})$.
\end{lemma}

\it{Proof.}\rm\quad
We will only give a sketch of the proof, which is similar to that of Lemma 3.13 of \cite{Morales-Rodrigo}.
For $k\geq \max\{ 2,\|p_0\|_{L^\infty(\Omega)}\}$, let  $q_k=\max\{q-k,0\}$ and $\Omega_k(t)=\{x\in \Omega: q(x,t)>k \}$.
Multiplying the equation of $q$ by $q_k$,  we obtain
\begin{equation}\label{2.10}
\begin{array}{ll}
&\displaystyle \frac{d}{dt}\int_{\Omega}q_k^{2}(c+1)^\alpha e^{\rho w}+2\int_{\Omega} |\nabla q_k |^2+2\int_{\Omega} q_k ^2+8\int_{\Omega}q_k^{2}(c+1)^\alpha e^{\rho w}\\[3mm]
\leq & \displaystyle c_1\int_{\Omega} q_k ^3+  c_1k\displaystyle \int_{\Omega} q_k^2 +c_1k^2 \int_{\Omega} q_k+
c_1 \displaystyle\int_{\Omega}(q_k^2+kq_k)  |\Delta c |
\end{array}
\end{equation}
for some $c_1>0$ independent of $k$.
By the boundedness of $q$ in $L^n(\Omega)$ for any $n>1$, the  Gagliardo--Nirenberg inequality and  Young inequality, we obtain
$$\displaystyle c_1\|q_k\|^3_{L^3(\Omega)}\leq \frac 14\|q_k\|^2_{H^1(\Omega)} +c_2\|q_k\|_{L^1(\Omega)},
$$
$$
\displaystyle c_1k\|q_k\|^2_{L^2(\Omega)}\leq \frac 14\|q_k\|^2_{H^1(\Omega)} +c_2k^2\|q_k\|_{L^1(\Omega)},$$
$$c_1 \displaystyle\int_{\Omega}q_k^2 |\Delta c |\leq \frac 14\|q_k\|^2_{H^1(\Omega)} +c_2\|q_k\|^2_{L^2(\Omega)}\|\Delta c\|^2_{L^2(\Omega)},$$
$$c_1 k\displaystyle\int_{\Omega}q_k |\Delta c |\leq \frac 14\|q_k\|^2_{H^1(\Omega)} +c_2k^2(1+\|\Delta c\|^8_{L^8(\Omega)})|\Omega_k|^{\frac 32}.$$
Inserting the above estimates into \eqref{2.10}, we have
$$
\begin{array}{ll}
&\displaystyle \frac{d}{dt}\int_{\Omega}q_k^{2}(c+1)^\alpha e^{\rho w}+\int_{\Omega} |\nabla q_k |^2+\int_{\Omega} q_k ^2+8\int_{\Omega}q_k^{2}(c+1)^\alpha e^{\rho w}\\[2mm]
\leq & c_2\|q_k\|^2_{L^2(\Omega)}\|\Delta c\|^2_{L^2(\Omega)}+c_2k^2(1+\|\Delta c\|^8_{L^8(\Omega)})|\Omega_k|^{\frac 32}
+(c_1+2c_2)k^2\|q_k\|_{L^1(\Omega)}\\
\leq & c_2\|q_k\|^2_{L^2(\Omega)}\|\Delta c\|^2_{L^2(\Omega)}+\displaystyle\frac 12\|q_k\|^2_{H^1(\Omega)}+c_3k^4(1+\|\Delta c\|^8_{L^8(\Omega)})|\Omega_k|^{\frac 32}.
\end{array}
$$
On the other hand, according  to the  relation between distribution functions and $L^p$ integrals (see e.g. (2.6) of \cite{Reyes}), we can see that
$$
(r+1)\int^\infty_0 s^r|\Omega_s(t)|ds=\|q(t)\|^{r+1}_{L^{r+1}(\Omega)}.
$$
Hence taking
into account Lemma 2.3, we get
 $$(k-1)^{16}|\Omega_k(t)|< \int^k_{k-1} s^{16}|\Omega_s(t)|ds<\int^\infty_0 s^{16}|\Omega_s(t)|ds\leq \frac 1{17}\|q(\cdot,t)\|^{17}_{L^{17}(\Omega)}$$
 and thus
\begin{equation}\label{2.11}
\begin{array}{ll}
&\displaystyle \frac{d}{dt}\int_{\Omega}q_k^{2}(c+1)^\alpha e^{\rho w}+8\int_{\Omega}q_k^{2}(c+1)^\alpha e^{\rho w}\\[2mm]
\leq & c_2\|\Delta c\|^2_{L^2(\Omega)}\displaystyle\int_{\Omega}q_k^{2}(c+1)^\alpha e^{\rho w}+c_4(1+\|\Delta c\|^8_{L^8(\Omega)})|\Omega_k|^{\frac 54}.
\end{array}
\end{equation}
Therefore if $h(t)=8-c_2\|\Delta c\|^2_{L^2(\Omega)}$, then
$$
\int_{\Omega}q_k^{2}(c+1)^\alpha e^{\rho w}\leq c_4 e^{-\int^t_0 h(s)ds}\int ^t_0(1+\|\Delta c\|^8_{L^8(\Omega)})e^{\int^s_0 h(\sigma)d\sigma}
|\Omega_k(s)|^{\frac 54}ds.
$$
Furthermore, since $ e^{-\int^t_0 h(s)ds}=e^{-8t}  e^{c_2\int^t_0 \|\Delta c\|^2_{L^2(\Omega)}ds} \leq c_5 e^{-8t} $ by Lemma 2.2  and $e^{\int^s_0 h(\sigma)d\sigma}\leq e^{8s} $,  we get
$$
\begin{array}{rl}
\displaystyle\int_{\Omega}q_k^{2}\leq & c_6 \displaystyle\int ^t_0 e^{-8(t-s)}(1+\|\Delta c\|^8_{L^8(\Omega)})
|\Omega_k(s)|^{\frac 54}ds\\[2mm]
\leq & c_6 \displaystyle \int ^t_0 e^{-8(t-s)}(1+\|\Delta c\|^8_{L^8(\Omega)})ds\cdot\displaystyle\sup_{t\geq 0}|\Omega_k(t)|^{\frac 54}.
\end{array}
$$
To estimate the integral term in the  right-side of the above inquality, we  apply Lemma \ref{lemma24} with $r=8$ and Lemma \ref{lemma23} to get
 $$
\int ^t_0 e^{-8(t-s)}\|\Delta c\|^8_{L^8(\Omega)}ds\leq c_7
$$ and thus
$\int_{\Omega}q_k^{2}\leq c_8 ( \displaystyle\sup_{t\geq 0}|\Omega_k(t)| )^{\frac 54}
$.

On the other hand,
 $\int_{\Omega}q_k^{2}(t)\geq  \int_{\Omega_j(t)}q_k^{2}(t)\geq (j-k)^2|\Omega_j(t)|$ for $j>k$.
 Consequently
 $$
 (j-k)^2\displaystyle\sup_{t\geq 0}|\Omega_j(t)|\leq c_8 ( \displaystyle\sup_{t\geq 0}|\Omega_k(t)| )^{\frac 54}
  |\Omega_k(t)|.
  $$
According to Lemma B.1 of \cite{Kinderlehrer}, there exists $k_0<\infty$ such that $|\Omega_{k_0}(t)|=0$ for all $t\in (0, T_{max})$. Therefore $\|q(\cdot,t)\|_{L^\infty(\Omega)}\leq k_0$ for any $t\in (0, T_{max})$ and thereby the proof is complete.

\vskip 3mm


{\it Proof of  Theorem 1.1.}~ By the boundedness of $p$ in $L^\infty((0, T_{max}),L^\infty(\Omega))$ from Lemma \ref{lemma25} and  a bootstrap argument as in \cite{Morales-Rodrigo},
we can see that the global existence of classical solutions to (1.1) is  an immediate consequence of Lemma \ref{lemma21},  i.e., $ T_{\max} = \infty$.
Indeed, supposed that $T_{\max}<\infty$, then by Lemma 3.15 and Lemma 3.19 of \cite{Morales-Rodrigo}, we can see that for any $s>N+2$ and $t\leq T_{\max} $ $$\|c(\cdot,t)\|_{W^{1,s}(\Omega)}+\|w(\cdot,t)\|_{W^{1,s}(\Omega)}\leq C.
$$
Further by Lemma 3.20 of \cite{Morales-Rodrigo}, we have $\|p(\cdot,t)\|_{W^{1,s}(\Omega)}\leq C$ which contradicts \eqref{2.1} and thus implies that  $ T_{\max} = \infty$.
Moreover, since $\tau:=\min\{1,\frac {T_{\max}}6\}=1 $, there exists a constant $C>0$ independent of time $t$ such that  $\|p(\cdot,t)\|_{L^\infty(\Omega)}\leq C$ for all $t\geq 0$ by retracing the proofs of Lemma \ref{lemma23} and Lemma \ref{lemma25}.  This completes the proof of Theorem 1.1.

\section{Proof of Theorem 1.2}
In this section, on the basis of the $L^\infty$-bound  of $p$ provided by Theorem \ref{theorem1.1},  we shall
look at the asymptotic behavior of 
the solution $(p,c,w)$ of the problem (1.1).
\subsection{$L^r$-convergence of solutions in two dimensions }
When either $w_0 > 1$ or $ \|w_0-1\|_{L^\infty(\Omega)}<\delta$ for some  $\delta> 0$, the authors of \cite{Morales-Rodrigo}
removed the time dependence of the $L^\infty$-bound of $p$ (see Lemma 5.8 of \cite{Morales-Rodrigo}) and
thereby investigated the  asymptotic behavior of  solutions to (1.1).
 In this subsection, 
 on the basis of the $L^\infty$-bound of $p$ being independent of time as provided by Theorem 1.1, we shall derive the same estimates as in Lemma  5.6 and Lemma 5.7 of \cite{Morales-Rodrigo} under the weaker assumption that $w_0>1-\frac 1 \rho$. We shall show that the solution $(p,c,w)$ to \eqref{1.1} converges to the  homogeneous  steady state $(1,0,1)$ as $t\rightarrow\infty$. 

Before going into the details, let us first collect some useful related estimates. It  should  be noted that no other assumptions on the initial data $(p_0, c_0, w_0)$ are made except for reasonable regularity, i.e., \eqref{1.4}.

\begin{lemma}\label{lemma32}\rm{(\cite[Lemmas 3.4, 5.1, 5.2, 3.8]{Morales-Rodrigo})}~\it Let $(p,c,w)$ be the global, classical solution of \eqref{1.1}. Then
\begin{equation}\label{3.1}
\begin{array}{l}
\left|\displaystyle\int_ 0^\infty \int_{\Omega} p(1-p)\right|\leq \max\{\int_\Omega p_0,|\Omega|\}/\lambda;
\end{array}
\end{equation}
\begin{equation}\label{3.2}
\displaystyle\int_ 0^\infty \int_{\Omega} p|w-1|
\leq \|w_0-1\|_{L^1(\Omega)};
\end{equation}
\begin{equation}\label{3.3}
\displaystyle\int_ 0^\infty \int_{\Omega} p|\nabla c|^2
< \infty 
\end{equation}
\begin{equation}
\displaystyle \int_{\Omega} |\nabla c(t)|^2\leq e^{-2t}\left(\int_{\Omega} |\nabla c_0|^2+\mu^2\|c_0\|^2_{L^\infty(\Omega)}\max\{\int_\Omega p_0,|\Omega|\}(t+\frac 1 \lambda)\right).
\end{equation}
\end{lemma}


\begin{lemma}\label{lemma33a}\rm Under the assumptions of Theorem 1.1, we have
\begin{equation}\label{3.5a}
\displaystyle \sup_{t\geq 0}\| c(t)\|_{W^{1,\infty}(\Omega)}
\leq C.
\end{equation}
\end{lemma}
\it{Proof.}\rm\quad
We know that $c$ solves the linear equation
$$
 c_t=\Delta c-c+f
$$
under the Neumann boundary condition with $f := -\mu pc $.  Since $p\geq 0$,  we know that $0\leq c(x,t)\leq$  $\|c_0\|_{L^\infty(\Omega)} e^{-t}$ by the standard sub-super solutions method. On the other hand, by Theorem \ref{theorem1.1}, $\displaystyle \sup_{t\geq 0}$ $\|p(t)\|_{L^\infty(\Omega)} $$ \leq c_1$, which readily implies that $\displaystyle
\|f\|_{L^\infty((0, \infty);L^\infty(\Omega))}\leq c_1$. Now upon a standard regularity argument we  can deduce the desired result. For the reader's convenience, we only give a brief sketch of the main ideas, and would like refer to the proof of Lemma 1 in \cite{Kowalczyk} or Lemma 4.1 in \cite{Horstmann} for more details. Indeed, according to the variation-of-constants formula of  $c$, we have for $t>2$
$$
c(\cdot,t)= e^{(t-1)(\Delta-1)}c(\cdot,1)+\int_{1}^t e^{(t-s)(\Delta-1)}f(\cdot,s)ds.
$$
So by Lemma \ref{lemma316}(ii), we infer that
\begin{equation*}
\begin{array}{rl}
\|\nabla c(\cdot,t)\|_{L^\infty(\Omega)}&\displaystyle\leq 2k_2 \|c(\cdot,1)\|_{L^1(\Omega)}+k_2\int_{1}^t (1+(t-s)^{-\frac 12}) e^{-(t-s)(\lambda_1+1)}\|f(\cdot,s)\|_{L^\infty(\Omega)}ds\\
& \leq  2k_2 \|c(\cdot,1)\|_{L^1(\Omega)}+k_2c_1\displaystyle\int_{0}^\infty (1+\sigma^{-\frac 12}) e^{-\sigma}d\sigma.
\end{array}
\end{equation*}

\vskip3mm

\begin{lemma}\label{lemma35}\rm{(\cite [Lemma 5.4] {Morales-Rodrigo})}~\it Let $(p,c,w)$ be the global, classical solution of \eqref{1.1}. Then for every $t\geq 0$ and $\kappa>0$,
\begin{equation}\label{3.5}
\displaystyle\frac d {dt} F(p(t), w(t))=G(p(t), w(t), c(t)),
\end{equation}
where
 $$
 F(p,w)=\kappa \int_\Omega |\nabla w|^2+ \int_\Omega p(\ln p-1) +\int_\Omega p(w-1)- \gamma \kappa  \int_\Omega p(w-1)^2,
$$
and
$$
\begin{array}{rl}
G(p, w, c)=&-\displaystyle \int_\Omega \displaystyle\frac {|\nabla p|^2}p+\int_\Omega \displaystyle\frac \alpha {1+c}\nabla p\cdot\nabla c+
\int_\Omega (2\alpha \gamma \kappa (1-w)+\alpha \rho) \frac p{1+c}\nabla c\cdot\nabla w\\[3mm]
&
+\displaystyle\int_\Omega (\rho^2-2\gamma \kappa +2\rho \gamma \kappa (1-w)) p |\nabla w|^2+ \lambda \int_\Omega p(1-p)\ln p\\
&+\displaystyle
\lambda \rho \int_\Omega p(1-p)(w-1) + \gamma \rho  \int_\Omega  p^2(1-w)\\[2mm]
&+\displaystyle
 2\gamma^2 \kappa \int_\Omega  p^2(w-1)^2- \lambda\gamma \kappa  \int_\Omega p(1-p)(w-1)^2.
\end{array}
$$
\end{lemma}

\begin{lemma}\label{lemma36} If $w_0>1-\displaystyle\frac 1 \rho$, then there exists  $\kappa>0$ such that
\begin{equation}\label{3.6}
G(p, w, c)\leq -\frac 12 \displaystyle \int_\Omega \displaystyle\frac {|\nabla p|^2}p-
\frac 12 \displaystyle \int_\Omega  p|\nabla w|^2
+C\int_\Omega p|w-1|+C\int_{\Omega} p|\nabla c|^2
\end{equation}
for some $C>0$.
\end{lemma}
\it{Proof.}\rm\quad
By the H\"{o}lder and  Young inequalities we have
$$
\int_\Omega \displaystyle\frac \alpha {1+c}\nabla p\cdot\nabla c\leq \frac 12 \displaystyle \int_\Omega \displaystyle\frac {|\nabla p|^2}p
+\frac {\alpha^2} 2 \int_{\Omega} p|\nabla c|^2,
$$
and
$$
\int_\Omega (2\alpha \gamma \kappa (1-w)+\alpha \rho) \frac p{1+c}\nabla c\cdot\nabla w
\leq \frac 12 \displaystyle \int_\Omega  p|\nabla w|^2
+c_1\int_{\Omega} p|\nabla c|^2
$$
for some $c_1>0$.

As $w_0>1-\displaystyle\frac 1 \rho$, we can find some $\varepsilon_1>0$ such that  $\rho(1-w_0)_+\leq 1-\varepsilon_1$, where $(1-w_0)_+=\max\{0,1-w_0\}$. Hence from the $w$-equation in \eqref{1.1}, it follows that
\begin{equation}\label{3.7}
1-w=(1-w_0)e^{-\gamma\int^t_0 p(s)ds},
\end{equation}
and thus
$$\begin{array}{rl}
\displaystyle\int_\Omega (\rho^2-2\gamma \kappa +2\rho \gamma \kappa (1-w)) p |\nabla w|^2
&\leq
\displaystyle\int_\Omega (\rho^2-2\gamma \kappa +2\rho \gamma \kappa (1-w_0)_+) p |\nabla w|^2 \\[3mm]
&\leq\displaystyle\int_\Omega (\rho^2-2\gamma \kappa \varepsilon_1) p |\nabla w|^2 \\
&\leq-\displaystyle\int_\Omega p |\nabla w|^2
\end{array}
$$
if we pick $\kappa>0$ sufficiently large such that $\rho^2-2\gamma \kappa \varepsilon_1<-1$.

Denote the lower-order terms of $G(p, w, c)$ by $\theta(p, w)$, i.e.,
$$
\begin{array}{rl}
\theta(p, w)=:&
\lambda\displaystyle \int_\Omega p(1-p)\ln p+\displaystyle
\lambda \rho \int_\Omega p(1-p)(w-1) + \gamma \rho  \int_\Omega  p^2(1-w)\\[2mm]
&+\displaystyle
 2\gamma^2 \kappa \int_\Omega  p^2(w-1)^2- \lambda\gamma \kappa  \int_\Omega p(1-p)(w-1)^2.
\end{array}
$$
Since $s(1 - s) \ln s \leq  0$ for $s\geq 0$,  we get
$$
\begin{array}{rl}
\theta(p, w)
\leq &
\displaystyle
\lambda \rho \int_\Omega p(1-p)(w-1) + \gamma \rho  \int_\Omega  p^2(1-w)\\[2mm]
&+\displaystyle
 2\gamma^2 \kappa \int_\Omega  p^2(w-1)^2- \lambda\gamma \kappa  \int_\Omega p(1-p)(w-1)^2\\[2mm]
 \leq &
 c_2(\|p\|_{L^\infty(\Omega)},\|w-1\|_{L^\infty(\Omega)})\displaystyle\int_\Omega  p|w-1|,\\[2mm]
\end{array}
$$
which, along with $\|p(\cdot,t)\|_{L^\infty(\Omega)}\leq C$ from Theorem \ref{theorem1.1} and $\|w(\cdot,t)-1\|_{L^\infty(\Omega)}\leq \|w_0-1\|_{L^\infty(\Omega)}$ from \eqref{3.7},   yields
$$ \theta(p, w)\leq  c_3 \displaystyle\int_\Omega  p|w-1|.
$$
The desired result \eqref{3.6} then immediately follows.

\vskip3mm

\begin{lemma}\label{lemma37} If  $w_0>1-\displaystyle\frac 1 \rho$, then
\begin{equation}\label{3.8}
\displaystyle\sup_{t\geq 0}\int_\Omega |\nabla w(t)|^2+
 \displaystyle \int_0^\infty\int_\Omega \displaystyle\frac {|\nabla p|^2}p+\displaystyle \int^\infty_0\int_\Omega  p|\nabla w|^2<\infty.
\end{equation}
\end{lemma}
\it{Proof.}\rm\quad
Combining Lemmas \ref{lemma35} and \ref{lemma36}, we have
\begin{equation}\label{3.9}
\displaystyle\frac d {dt} F(p(t), w(t))+\frac 12 \displaystyle \int_\Omega \displaystyle\frac {|\nabla p|^2}p+
\frac 12 \displaystyle \int_\Omega  p|\nabla w|^2\leq
C\int_\Omega p|w-1|+C\int_{\Omega} p|\nabla c|^2.
\end{equation}
Hence \eqref{3.8} follows upon integration on the time variable, 
and using \eqref{3.2} and \eqref{3.3}.

\vskip3mm

\begin{lemma}\label{lemma38} If $w_0>1-\displaystyle\frac 1 \rho$, then for any $r\geq 2$
\begin{equation}\label{3.10}
\displaystyle\lim_{t\rightarrow \infty}\|p(\cdot,t)-\overline p(t)\|_{L^r(\Omega)}=0,
\end{equation}
\begin{equation}\label{3.11}
\displaystyle\lim_{t\rightarrow \infty}|\overline p(t)-1|=0,
\end{equation}
where $\overline p(t)=\frac 1 {|\Omega|}\int_\Omega p(\cdot, t)$, and
\begin{equation}\label{3.12}
\displaystyle\lim_{t\rightarrow \infty}\|w(\cdot,t)-1\|_{L^r(\Omega)}=0.
\end{equation}
\end{lemma}

\it{Proof.}\rm\quad
The proofs of \eqref{3.10} and \eqref{3.11} are similar to those of Lemma 5.9--5.11 of \cite{Morales-Rodrigo} respectively. However, for the reader's convenience, we only give a brief sketch of \eqref{3.11}.
In fact, from \eqref{1.1} and the Poincar\'{e}--Wirtinger inequality, it follows that 
\begin{equation*}\label{4.12}
\begin{array}{rl}
\overline p_t= & \lambda (\displaystyle\overline p - \overline p^2- \displaystyle\frac {1}{|\Omega|}\int_\Omega (p-\overline p)^2\\[2mm]
\geq & \lambda \overline p(\displaystyle 1-\overline p-c_1\int_\Omega  \frac {|\nabla p|^2}p).
\end{array}
\end{equation*}
Hence by \eqref{3.8}, we get \begin{equation*}
\begin{array}{rl}
\overline p(t)\geq  & \displaystyle\overline p_0 \exp \{\lambda t -\lambda\int ^t_0 \overline p(s)ds-c_1\lambda \int ^\infty_0\int_\Omega  \frac {|\nabla p|^2}p ds \} \\[2mm]
\geq &  \displaystyle c_2 \exp \{\lambda t -\lambda\int ^t_0 \overline p(s)ds \},
\end{array}\end{equation*}
which means that \eqref{3.11} is valid due to either $\overline p(t)\rightarrow 1$ or $\overline p(t)\rightarrow 0$ in Lemma 5.10 of \cite{Morales-Rodrigo}.
Indeed, supposed that $\overline p(t)\rightarrow 0$, then there exists $t_0>1$ such that
$
\overline p(t)\leq \frac 12
$  and thus $ \int ^t_0 \overline p(s)ds \leq \frac t 2+ \int ^{t_0}_0 \overline p(s)ds$
for all $t\geq t_0$. 
Therefore we arrive at  $\overline p(t)\geq c_3 e^\frac{\lambda t}2$ for all $t\geq t_0$, which contradicts  $\overline p(t)\rightarrow 0$.

Now we turn to show \eqref{3.12}.
Invoking the Poincar\'{e} inequality in the form
$$\displaystyle\int_\Omega |\varphi(x)-\frac 1 {|\Omega|}\int_\Omega \varphi(y)dy|^2 dx \leq C_p \int_\Omega | \nabla \varphi|^2 dx~~\hbox{
for all~} \varphi\in W^{1,2}(\Omega)$$
 for some $C_p > 0$, one can find that for all $j\in \mathbb{N}$
\begin{equation}\label{3.13}
  \begin{array}{rl}
  \displaystyle\int^{j+1}_j \|p(s)- \overline p(s)\|^2_{L^2(\Omega)}ds &\leq C_p\displaystyle
\int^{j+1}_j  \|\nabla p(s)\|^2_{L^2(\Omega)}ds\\[2mm]
& \leq C_p\displaystyle\sup_{t\geq 0}\|p(t)\|_{L^\infty(\Omega)}
\displaystyle \int^{j+1}_j  \int_\Omega \frac{|\nabla p(s)|^2}{p(s)}  ds, \\
\end{array}
\end{equation}
which,  along with \eqref{3.8} and Theorem 1.1, shows that
\begin{equation}\label{3.14}
\int_\Omega\int^{j+1}_j |p(x,s)- \overline p(s)|^2 dsdx =\int^{j+1}_j \|p(s)- \overline p(s)\|^2_{L^2(\Omega)}ds\rightarrow 0 \end{equation}
as $j\rightarrow\infty$.

Now defining $p_j(x):=\int^{j+1}_j |p(x,s)- \overline p(s)|^2 ds$, $x\in \Omega, j\in \mathbb{N}$,  \eqref{3.14} tells us that
$p_j \rightarrow 0$ in $L^1(\Omega)$ as $j\rightarrow\infty$.
There exist a certain null set $Q\subseteq\Omega$ 
 and a subsequence $(j_k)_{k\in \mathbb{N}}\subset \mathbb{N}$ such that $j_k\rightarrow\infty$ and
$p_{j_k}(x) \rightarrow 0$ for every $x\in \Omega\setminus Q$ 
as $k\rightarrow\infty$.
Restated in the original variable, this becomes
\begin{equation}\label{3.15}
\int^{j_k+1}_{j_k} |p(x,s)- \overline p(s)|^2 ds \rightarrow 0
\end{equation}
for every $x\in \Omega\setminus Q$
as $k\rightarrow\infty$.

Therefore, from \eqref{3.7} and $p(x,t)\geq 0$, it follows that for any $x\in \Omega\setminus Q$
\begin{equation}\label{3.16}
\begin{array}{rl}
|w(x,t)-1|&
\leq \|w_0-1\|_{L^\infty(\Omega)}\exp\{-\gamma\int^{[t]}_0 p(x,s)ds\}\\[2mm]
&
\leq \|w_0-1\|_{L^\infty(\Omega)}\exp\{-\gamma\sum_{k=0}^{m(t)}\int^{j_k+1}_{j_k}p(x,s)ds\}\\[2mm]
&\leq \|w_0-1\|_{L^\infty(\Omega)}\exp\{\gamma\sum_{k=0}^{m(t)}\int^{j_k+1}_{j_k}|p(x,s)- \overline p(s)|ds-
\gamma\sum_{k=0}^{m(t)}\int^{j_k+1}_{j_k} \overline p(s)ds\}\\[2mm]
&\leq \|w_0-1\|_{L^\infty(\Omega)}\exp\{\gamma\sum_{k=0}^{m(t)}(\int^{j_k+1}_{j_k}|p(x,s)- \overline p(s)|^2 ds)^{\frac12}\\[2mm]
& \qquad -\gamma\sum_{k=0}^{m(t)}\int^{j_k+1}_{j_k} \overline p(s)ds\},
\end{array}
\end{equation}
where $m(t):=\displaystyle\max_{k\in \mathbb{N}}\{j_k+1,[t]\}$.
Furthermore,
by \eqref{3.11}, there exists $k_0\in \mathbb{N}$ such that $\int^{j_k+1}_{j_k} \overline p(s)ds\geq \frac 12 $ for all $k\geq k_0$.
Hence by the fact that $m(t)\rightarrow \infty$ as $t\rightarrow \infty$ and \eqref{3.15}, we  obtain that  $w(x,t)-1\rightarrow 0$  almost
everywhere in $\Omega$ as $t\rightarrow \infty$. On the other hand,
as $|w(x,t)-1|\leq \|w_0-1\|_{L^\infty(\Omega)}$, the dominated convergence theorem ensures that \eqref{3.12} holds for any $r\in (2,\infty)$.

\vskip 3mm

{\bf Remark 3.1}.  1) It is observed that since $W^{1,2}(\Omega)\hookrightarrow L^\infty(\Omega)$ is invalid in the two-dimensional setting, $\|p(s)- \overline p(s)\|^2_{L^2(\Omega)}$  in  \eqref{3.13} cannot be replaced by $\|p(s)- \overline p(s)\|^2_{L^\infty(\Omega)}$, and thus we cannot infer that
 $\displaystyle\lim_{t\rightarrow \infty}\|w(\cdot,t)-1\|_{L^\infty(\Omega)}=0$, even though we have established that all the related estimates of $(p,c,w)$
  in \cite{Morales-Rodrigo} continue to hold under the milder condition imposed on the initial data $w_0$.

2) Similarly to the remark above, we note that, even though $\|w(\cdot,t)\|_{W^{1,n}(\Omega)}\leq C(T)$ for any $n\geq 2$ and $t\leq T$,
we are not able to infer the global estimate $\sup_{t\geq 0}\int_\Omega|\nabla w(t)|^{2+\varepsilon}$ $\leq C$. Otherwise, we would be able to apply regularity estimates
for bounded solutions of semilinear parabolic equations (see \cite{Porzio} for instance) to obtain the H\"{o}lder estimates of $p(x,t)$ in $\Omega\times (1,\infty)$, and thereby conclude
$\displaystyle\lim_{t\rightarrow \infty}\|p(\cdot,t)-1\|_{L^\infty(\Omega)}=0$.
As things stand at the moment, we are only able to infer convergence in $L^r$.

\subsection{$L^\infty$-convergence of solutions with exponential rate in one dimension}

It is observed that the results, in particular Lemma \ref{lemma38}, in the previous subsection are still valid in the one-dimensional case. Moreover, in the one-dimensional setting, the weak convergence result  in Lemma \ref{lemma38} can be improved via a bootstrap argument. In fact,  we shall
derive some a priori estimates of $(p,c,w)$ and thereby demonstrate that
$(p,c,w)$ converges to $(1,0,1)$ in $L^\infty(\Omega)$ as $t\rightarrow \infty$. Furthermore, by a regularity argument involving the variation-of-constants formula
for $p$ and smoothing $L^p-L^q$ type estimates for the Neumann heat semigroup, we will show that $p(\cdot,t)-1$ decays exponentially in $L^\infty(\Omega)$.

As pointed out in the Introduction, the main technical difficulty in the derivation of Theorem 1.2 stems from the coupling between $p$
and $w$. Indeed, the lack of regularization effect in the space variable in the $w$-equation
and the presence of $p$ there demand tedious estimates of the solution.

The following lemma  plays a crucial role in establishing the uniform convergence of $p$ as $t\rightarrow \infty$  (see Lemma \ref{lemma312}).
Thought the proof thereof only  involves elementary analysis, we give a full proof here for the sake of the reader's convenience since we could not find a precise reference covering our situation. 
\begin{lemma}\label{lemma39}  Let $k(t)$ be a function satisfying
$$
k(t)\geq 0, \quad \int^\infty_0k(t)dt<\infty.
$$
If $k'(t)\leq h(t)$ for some $h(t)\in L^1(0,\infty)$, then $k(t)\rightarrow 0$ as $t\rightarrow \infty$.
\end{lemma}
\it{Proof.}\rm\quad
Supposing the contrary, then we can find $A>0$ and a sequence
$(t_j)_{j\in \mathbb{N}}\subset (1,\infty)$ such that  $t_j\geq t_{j-1}+2$, $t_j \rightarrow \infty$  as $j\rightarrow \infty$ and
$k(t_j)\geq A$ for all $j\in \mathbb{N} $. On the other hand, by $k'(t)\leq h(t)$, we have
\begin{equation}\label{3.17}
k(t_j-\tau)\geq k(t_j)-\int^{t_j}_{t_j-\tau}|h(s)|ds\geq k(t_j)-\int^{t_j}_{t_j-1}|h(s)|ds
\end{equation}
for all $\tau \in (0,1)$.

Since $h(t)\in L^1(0,\infty)$,  we have $\int^{t_j}_{t_j-1}|h(s)|ds\rightarrow 0 $ as $t_j \rightarrow \infty$ and thereby there exists $j_0\in \mathbb{N} $ such that
 $\int^{t_j}_{t_j-1}|h(s)|ds \leq \frac A2 $
 for all $j\geq j_0$, which along with \eqref{3.17} implies that
\begin{equation}\label{3.18}
k(t_j-\tau)\geq  k(t_j)-\frac A2 \geq \frac A2
\end{equation}
for $j\geq j_0$ and $\tau \in (0,1)$.  It follows that $\int^{t_j}_{t_j-1}k(t)dt\geq \frac A2$ for all $j\geq j_0$, which contradicts $\int^{\infty}_{0}k(t)dt< \infty$ and thus completes the proof of the lemma.

\vskip 3mm

\begin{lemma}\label{lemma310} If  $w_0>1-\displaystyle\frac 1 \rho$, then there exists a constant  $C>0$ such that
\begin{equation}\label{3.19}
 \displaystyle \int_0^\infty\int_\Omega \displaystyle e^{\rho w} {| z_x|^2}\leq C
\end{equation}
where $z=p e^{-\rho w} $.
\end{lemma}
\it{Proof.}\rm\quad  We know that $ z_x=e^{-\rho w}  p_x-\rho z w_x$ and thus
$$
e^{\rho w}|z_x|^2\leq 4e^{-\rho w} |p_x|^2+4p^2 e^{-\rho w} \rho^2 | w_x|^2.
$$
Integrating over $\Omega\times (0,\infty)$ and taking \eqref{2.5a} into account, we have
$$
\begin{array}{rl}
\displaystyle \int_0^\infty\int_\Omega \displaystyle e^{\rho w} {| z_x|^2} &\leq
 4\displaystyle \int_0^\infty\int_\Omega \displaystyle | p_x|^2+4\rho^2 \displaystyle\sup_{t\geq 0}\|p(t)\|_{L^\infty(\Omega)}  \displaystyle \int_0^\infty\int_\Omega p | w_x|^2\\[3mm]
 &\leq 4(1+ \rho^2)\displaystyle\sup_{t\geq 0}\|p(t)\|_{L^\infty(\Omega)}
 (
 \displaystyle \int_0^\infty\int_\Omega \displaystyle\frac {| p_x|^2}p+\displaystyle \int^\infty_0\int_\Omega  p| w_x|^2
 ).
\end{array}
$$
Hence by  Theorem \ref{theorem1.1} and Lemma \ref{lemma37}, we  get \eqref{3.19}.

\vskip 3mm

\begin{lemma}\label{lemma311} If  $w_0>1-\displaystyle\frac 1 \rho$, then there exists a constant  $C>0$ such that
\begin{equation}\label{3.20}
\begin{array}{rl}
&\displaystyle\frac d{dt}

\int_\Omega \displaystyle e^{\rho w} {| z_x|^2}+ \frac 13\displaystyle \int_\Omega e^{\rho w}z_t^2 \\[3mm]
 \leq & C (
\displaystyle \int_\Omega e^{\rho w} {| z_x|^2}+\int_\Omega  p| w_x|^2+\int_\Omega|c_{xx}|^2+\int_\Omega  p| c_x|^2+
 \int_\Omega  p|w-1|+ \int_\Omega  p(p-1)^2
 )
  \end{array}
\end{equation}
with $z=p e^{-\rho w} $.
\end{lemma}
\it{Proof.}\rm\quad  Note that $z$ satisfies
$$
z_t = e^{-\rho w}(e^{\rho w } z_x)_x -e^{-\rho w} (\displaystyle \frac{ z_x e^{\rho w}}{1+c}
\nabla c)_x
+ \lambda z (1-ze^{\rho w})-\rho\gamma e^{\rho w}z^2(1-w).
$$
Multiplying  the above equation by $z_t e^{\rho w} $ and integrating in the spatial
variable, we obtain
\begin{equation}
\label{3.21}
\begin{array}{rl}
& \displaystyle \int_\Omega e^{\rho w}z_t^2+\int_\Omega e^{\rho w } z_x z_{xt} \\[2mm]
=&- \displaystyle \int_\Omega e^{\rho w}z_t(\displaystyle \frac{\alpha}{1+c} z_x c_x +
\displaystyle \frac{\alpha z \rho}{1+c} w_x
 c_x-
\displaystyle \frac{\alpha z }{(1+c)^2}
| c_x|^2
+\displaystyle \frac{\alpha z }{1+c}
 c_{xx})\\[2mm]
&
+\displaystyle \int_\Omega e^{\rho w}z_t(\lambda z (1-ze^{\rho w})-\rho\gamma e^{\rho w}z^2(1-w)).
\end{array}
\end{equation}
Notice that

$$
\begin{array}{rl}
\displaystyle\int_\Omega e^{\rho w } z_x z_{xt}=&
\displaystyle\frac 12\frac d{dt}
 \int_\Omega \displaystyle e^{\rho w} {| z_x|^2} - \displaystyle\frac {\gamma \rho}2 \int_\Omega e^{2\rho w}z(1-w)| z_x|^2\\[3mm]
\geq &\displaystyle\frac 12\frac d{dt}
 \int_\Omega \displaystyle e^{\rho w} {| z_x|^2} -\displaystyle\frac {\gamma \rho}2\displaystyle\sup_{t\geq 0}\|p(t)\|_{L^\infty(\Omega)} \|1-w_0\|_{L^\infty(\Omega)}\displaystyle\int_\Omega \displaystyle e^{\rho w} {| z_x|^2},
\end{array}
$$
$$
- \displaystyle \int_\Omega z_t\displaystyle \frac{\alpha e^{\rho w}}{1+c} z_x c_x
\ \leq \ \frac 16\displaystyle \int_\Omega e^{\rho w}z_t^2
+c_1\displaystyle\sup_{t\geq 0}\| c_x (t)\|^2_{L^\infty(\Omega)} \displaystyle \int_\Omega e^{\rho w}| z_x |^2,
$$
$$
\displaystyle \int_\Omega z_t\displaystyle\frac{\alpha z  e^{\rho w}\rho}{1+c} w_x c_x
\ \leq \ \frac 16\displaystyle \int_\Omega e^{\rho w}z_t^2
+c_1\displaystyle\sup_{t\geq 0}\|  c_x (t)\|^2_{L^\infty(\Omega)} \displaystyle\sup_{t\geq 0}\|p(t)\|_{L^\infty(\Omega)}\displaystyle \int_\Omega p| w_x |^2,
$$
$$
\displaystyle \int_\Omega \displaystyle z_t
\displaystyle \frac{\alpha z e^{\rho w}}{(1+c)^2}
| c_x |^2
\ \leq \ \frac 16\displaystyle \int_\Omega e^{\rho w}z_t^2
+c_1\displaystyle\sup_{t\geq 0}\|p(t)\|_{L^\infty(\Omega)}\displaystyle\sup_{t\geq 0}\| c_x (t)\|^2_{L^\infty(\Omega)}
\displaystyle\int_\Omega p| c_x |^2,
$$
$$
-\displaystyle \int_\Omega \displaystyle z_t \displaystyle \frac{\alpha ze^{\rho w} }{1+c}
 c_{xx}  \ \leq \ \frac 16\displaystyle \int_\Omega e^{\rho w}z_t^2
+c_1\displaystyle\sup_{t\geq 0}\|p(t)\|^2_{L^\infty(\Omega)}\int_\Omega| c_{xx} |^2,
$$
$$
\begin{array}{rl}
 &\displaystyle \int_\Omega e^{\rho w}z_t(\lambda z (1-ze^{\rho w})-\rho\gamma e^{\rho w}z^2(1-w))\\[3mm]
=&\lambda \displaystyle \int_\Omega pz_t (1-p)-\rho\gamma\int_\Omega \displaystyle z_tp^2(1-w))\\[3mm]
\leq &\displaystyle \frac 16\displaystyle \int_\Omega e^{\rho w}z_t^2
+c_1  \displaystyle \sup_{t\geq 0}\|p(t)\|_{L^\infty(\Omega)}\int_\Omega p(1-p)^2+c_1 \displaystyle\sup_{t\geq 0}\|p(t)\|^3_{L^\infty(\Omega)} \|1-w_0\|_{L^\infty(\Omega)} \int_\Omega p|1-w|.
\end{array}
$$
Applying Theorem \ref{theorem1.1}, \eqref{3.5a} and
inserting the above inequalities  into \eqref{3.21}, we obtain \eqref{3.20}.

\vskip 3mm

Now we focus our attention on the decay properties of the solutions. Indeed, we will show that $p(x,t)$ converges  to $1$   with respect to the norm in $L^\infty(\Omega)$ as $t\rightarrow \infty$. Subsequently, we will establish 
the exponential decay of  $\|p(\cdot,t)-1\|_{L^\infty(\Omega)}$ with explicit rate.

\begin{lemma}\label{lemma313} If $w_0>1-\displaystyle\frac 1 \rho$, then
\begin{equation}\label{3.22a}
\displaystyle\lim_{t\rightarrow \infty}\|w(\cdot,t)-1\|_{L^\infty(\Omega)}=0.
\end{equation}
\end{lemma}

\it{Proof.}\rm\quad
From \eqref{3.7}, it follows that for any $\epsilon>0$
\begin{equation}\label{3.23}
\begin{array}{rl}
|w(x,t)-1|&\leq \|w_0-1\|_{L^\infty(\Omega)}\exp\{-\gamma\int^t_0 p(s)ds\}\\[2mm]
&\leq \|w_0-1\|_{L^\infty(\Omega)}\exp\{\gamma\int^t_0\|p(s)- \overline p(s)\|_{L^\infty(\Omega)}ds-\gamma\int^t_0 \overline p(s)ds\}\\[2mm]
&\leq \|w_0-1\|_{L^\infty(\Omega)}\exp\{\frac \gamma \epsilon\int^t_0 \|p(s)- \overline p(s)\|^2_{L^\infty(\Omega)}ds+\epsilon \gamma t-\gamma\int^t_0 \overline p(s)ds\},
\end{array}
\end{equation}
where $\overline p(t)=\frac 1 {|\Omega|}\int_\Omega p(\cdot, t)$.

On the other hand, by the Poincar\'{e}--Wirtinger inequality, the Sobolev imbedding theorem in one dimension and \eqref{3.8}, we have
\begin{equation}\label{3.24}
  \begin{array}{rl}
  \displaystyle\int^t_0 \|p(s)- \overline p(s)\|^2_{L^\infty(\Omega)}ds &\leq c_1\displaystyle
\int^\infty_0 \| p_x (s)\|^2_{L^2(\Omega)}ds\\[2mm]
& \leq c_1\displaystyle\sup_{t\geq 0}\|p(t)\|_{L^\infty(\Omega)}
\displaystyle\int^\infty_0 \int_\Omega \frac{| p_x (s)|^2}{p(s)}  ds\\
& \leq c_2
\end{array}
\end{equation}
for some constant $c_2>0$.
Combining \eqref{3.23} with \eqref{3.24} yields
\begin{equation}\label{3.25}
\|w(t)-1\|_{L^\infty(\Omega)}\leq  \|w_0-1\|_{L^\infty(\Omega)}\exp\{\frac {c_2 \gamma} \epsilon +\epsilon  \gamma t-\gamma \int^t_0 \overline p(s)ds\}
\end{equation}for $t\geq 0$. 
The assertion now follows from the last inequality and the proof is complete.

\vskip 3mm


\begin{lemma}\label{lemma312} If  $w_0>1-\displaystyle\frac 1 \rho$, then
\begin{equation}\label{3.26}
\displaystyle\lim_{t\rightarrow \infty}
\|p(\cdot,t)-1\|_{L^\infty(\Omega)}=0.
\end{equation}
\end{lemma}
\it{Proof.}\rm\quad
We first show that
\begin{equation}\label{3.27}
\displaystyle\lim_{t\rightarrow \infty}\|z(\cdot,t)-\overline{z}(t)\|_{L^\infty(\Omega)}=0
\end{equation}
where $\overline z(t)=\frac 1 {|\Omega|}\int_\Omega z(\cdot, t)$.
To this end, we consider the function $k(t)\geq 0$ defined by
$
k(t)=\int_\Omega \displaystyle e^{\rho w} {| z_x |^2}
$
and prove that
 \begin{equation}\label{3.28}
\displaystyle\lim_{t\rightarrow \infty}k(t)=0.
\end{equation}
By Lemmas \ref{lemma39},  \ref{lemma310}  and   \ref{lemma311}, it is enough to prove that
$$
h(t):=\displaystyle \int_\Omega e^{\rho w} {| z_x |^2}+\int_\Omega  p| w_x |^2+\int_\Omega| c_{xx} |^2+\int_\Omega  p| c_x |^2+
 \int_\Omega  p|w-1|+ \int_\Omega  p(p-1)^2\in L^1(0,\infty).
$$
Noting \eqref{3.19},  \eqref{3.8},  \eqref{3.3},  \eqref{3.2} and  \eqref{2.4a},  it  remains to  estimate
$\int^\infty_0\int_\Omega  p(p-1)^2 $.
In fact, multiplying the $p$-equation in \eqref{1.1}   by $p- 1$, we have
\begin{equation}\label{3.29}
\begin{array}{rl}
\displaystyle\frac12\frac d{dt}
\int_\Omega \displaystyle (p-1)^2= & - \displaystyle \int_\Omega | p_x |^2+\rho\int_\Omega  p w_x p_x 
+\alpha\int_\Omega  \displaystyle\frac  p{1+c}c_x  p_x 
-\lambda \int_\Omega p(p-1)^2\\[3mm]
 \leq & -\displaystyle\frac12 \displaystyle \int_\Omega | p_x |^2+C
 ( \int_\Omega  p^2| w_x |^2+ \int_\Omega  p^2| c_x |^2)
 -\lambda \int_\Omega p(p-1)^2.
  \end{array}
\end{equation}
 Hence, by the boundedness of $p$,  
  \eqref{3.3} and  \eqref{3.8}, we easily infer that $\int^\infty_0\int_\Omega  p(p-1)^2\leq C $.

 Furthermore, by the Poincar\'{e}--Wirtinger inequality and the Sobolev imbedding theorem in one dimension, we have
\begin{equation}
\label{3.30}
   \displaystyle \|z(t)- \overline z(t)\|_{L^\infty(\Omega)} \leq C_p\displaystyle
 \| z_x (t)\|_{L^2(\Omega)},
\end{equation}
which along with   \eqref{3.28} yields  \eqref{3.27}.

On the other hand, for any $\{t_j\}_{j\in \mathbb{N}}\subset (1,\infty)$,   there exists a subsequence along which $z(\cdot,t_j)-e^{-\rho}\rightarrow 0$ a.e. in $ \Omega $ as $j\rightarrow \infty$ by Lemma 3.6. We apply the dominated convergence theorem along with the uniform majorization $|z(\cdot,t_j)|\leq \displaystyle\sup_{j \geq 1}\|z(t_j)\|_{L^\infty(\Omega)}\leq C$ to infer  that
\begin{equation}\label{3.31}
\displaystyle\lim_{t\rightarrow \infty}|\overline{z}(t)-e^{-\rho}|=0.
\end{equation}
Hence
 $$
 \begin{array}{rl}
 \|p(\cdot,t)-1\|_{L^\infty(\Omega)}= &\|e^{\rho w}z-1\|_{L^\infty(\Omega)}\\
 \leq & e^{\rho(1+\|w_0\|_{\infty(\Omega)})}(\|z(\cdot,t)-\overline z(t)\|_{L^\infty(\Omega)}+|\overline z(t)-e^{-\rho}|)+
 \|e^{\rho w(\cdot,t)}-e^{\rho}\|_{L^\infty(\Omega)}\\
 \leq & e^{\rho(1+\|w_0\|_{\infty(\Omega)})}(\|z(\cdot,t)-\overline z(t)\|_{L^\infty(\Omega)}+|\overline z(t)-e^{-\rho}|)+
c_1\| w(\cdot,t)-1\|_{L^\infty(\Omega)}
 \end{array}
 $$
 for some $c_1>0$, which, together with \eqref{3.22a}, \eqref{3.27} and  \eqref{3.31}, yields the desired result.

\vskip 3mm


Having established that
$p(x,t)$ converges to 1 uniformly with respect to $x\in \Omega$ as $t\rightarrow \infty$,
we now go on 
to establish an explicit exponential convergence rate.
Using \eqref{3.26}, we first look into the decay of $\int_\Omega | w_x (t)|^2$.

\begin{lemma}\label{lemma314} Let $w_0>1-\displaystyle\frac 1 \rho$. Then for any $\epsilon>0$, there exists $C(\epsilon)>0$ such that
\begin{equation}\label{3.32}
\displaystyle \int_\Omega | w_x (t)|^2 \leq C(\epsilon)e^{-2\gamma (1-\epsilon)t}.
\end{equation}
\end{lemma}
\it{Proof.}\rm\quad
From \eqref{3.7}, it follows that
$$
\begin{array}{rl}
 | w_x (t)|^2 \leq & 4 | w_{0x} |^2 e^{-2\gamma\int^t_0 p(s)ds}+4\gamma^2| w_0-1|^2 e^{-2\gamma\int^t_0 p (s)ds}
\left (\displaystyle\int^t_0 | p_x(s)|ds\right)^2\\[3mm]
\leq&  4 |w_{0x} |^2 e^{-2\gamma\int^t_0 p(s)ds}+4 t\gamma^2| w_0-1|^2 e^{-2\gamma\int^t_0 p(s)ds}
\displaystyle\int^t_0 | p_x(s)|^2 ds.
\end{array}
$$
Taking Lemma \ref{lemma312} into account, we know that for any $\epsilon>0$, there exists $t_\epsilon>1$ such that
$ p(x,t)>1-\epsilon $ for all $x\in \Omega, t> t_\epsilon$. Therefore
integrating the above inequality  in the space variable  yields
$$
\begin{array}{rl}
\displaystyle \int_\Omega | w_{x}(t)|^2 \leq & 4 e^{-2\gamma (1-\epsilon)(t-t_\epsilon)}\displaystyle \int_\Omega  |w_{0x}|^2+4 t\gamma^2\| w_0-1\|_{L^\infty(\Omega)}^2 e^{-2\gamma(1-\epsilon)(t-t_\epsilon)}
\displaystyle\int^\infty _0 \displaystyle \int_\Omega| p_{x}|^2\\[3mm]
\leq&  c_1(\epsilon)(1+t) e^{-2\gamma (1-\epsilon)t}
(\displaystyle\int_\Omega  |w_{0x}|^2  +
\| w_{0}-1\|_{L^\infty(\Omega)}^2 \displaystyle\sup_{t\geq 0}\| p(t)\|_{L^\infty(\Omega)}
\displaystyle\int^\infty _0 \displaystyle \int_\Omega\frac{| p_x|^2}{p}),

\end{array}
$$
for all $t>t_\epsilon$, which, along with \eqref{3.8},  implies \eqref{3.32}.

Now we utilize the decay  properties of  $ \int_\Omega | c_x(t)|^2$, $ \int_\Omega |w_x(t)|^2$  and the uniform convergence of
$ |p(x,t)-1| $  asserted by Lemma 3.11 to establish the decay property of $ \|p(\cdot,t)-1\|_{L^2(\Omega)} $.

\vskip 3mm

\begin{lemma}\label{lemma315} Let $w_0>1-\displaystyle\frac 1 \rho$. Then for any $\epsilon\in(0,\min\{1,\gamma,\lambda\})$, there exists $C(\epsilon)>0$ such that
\begin{equation}\label{3.33}
\displaystyle \|p(\cdot,t)-1\|_{L^2(\Omega)} \leq C(\epsilon)e^{-(\min\{1,\gamma,\lambda\}-\epsilon) t}.
\end{equation}
\end{lemma}
\it{Proof.}\rm\quad
By \eqref{3.26},  we know that for any $\epsilon\in(0,\min\{1,\gamma,\lambda\})$, there exists $t_\epsilon>1$ such that
$ p(x,t)>1-\epsilon $ for all $x\in \Omega, t> t_\epsilon$. Hence, we
 multiply the $p$-equation in \eqref{1.1}  by $p-1$ and integrate the result over $\Omega$ to get
  \begin{equation}\label{3.34}
  \begin{array}{rl}
\displaystyle\frac12\frac d{dt}
\int_\Omega \displaystyle (p-1)^2 &= - \displaystyle \int_\Omega | p_x|^2+\rho\int_\Omega  pw_x  p_x
+\alpha\int_\Omega  \displaystyle\frac  p{1+c} c_x p_x
-\lambda \int_\Omega p(p-1)^2\\[3mm]
 &\leq -\displaystyle\frac12 \displaystyle \int_\Omega | p_x|^2+c_1
 ( \int_\Omega  | w_x|^2+ \int_\Omega  | c_x|^2)
 -\lambda (1-\epsilon)\int_\Omega (p-1)^2
  \end{array}
\end{equation} for all $ t> t_\epsilon$.
Now, applying the Gronwall inequality, (3.4) and  Lemma \ref{lemma314}, we have
  $$
  \begin{array}{rl}
 \displaystyle\int_\Omega \displaystyle (p(t)-1)^2
 &
 \leq e^{-2\lambda(1-\epsilon)(t-t_\epsilon)}\displaystyle\int_\Omega \displaystyle (p(t_\epsilon)-1)^2
 +c_1\int^t_{0} e^{-2\lambda(1-\epsilon)(t-s)}
 (\displaystyle \int_\Omega  | w_x|^2+ \int_\Omega  | c_x|^2)\\[3mm]
 &
 \leq c_2(\epsilon)e^{-2\lambda(1-\epsilon)t}+c_3(\epsilon)\displaystyle\int^t_{0} e^{-2\lambda(1-\epsilon)(t-s)}
 (e^{-2\gamma(1-\epsilon)s}+ e^{-2(1-\epsilon)s}) \\
&
 \leq  c_4(\epsilon)e^{-2\min\{\lambda,1,\gamma\}(1-\epsilon)t},
 \end{array}
 $$
 where $c_i(\epsilon)>0$ $(i=2,3,4)$ are  independent of time $t$. This completes the proof.

\vskip 3mm
Moving forward, on the basis of Lemma \ref{lemma315}, we come to establish the exponential decay  of $ \displaystyle \|p(\cdot,t)-\overline p(t)\|_{L^\infty(\Omega)}$ by means of  a variation-of-constants representation of $p$, as follows:

\begin{lemma}\label{lemma317} Let $w_0>1-\displaystyle\frac 1 \rho$. Then for any $\epsilon\in(0,\min\{\lambda_1,1,\gamma,\lambda\})$, there exists $C(\epsilon)>0$ such that
\begin{equation}\label{3.35}
\displaystyle \|p(\cdot,t)-\overline p(t)\|_{L^\infty(\Omega)} \leq C(\epsilon)e^{-(\min\{\lambda_1,1,\gamma,\lambda\}-\epsilon) t}.
\end{equation}
\end{lemma}
\it{Proof.}\rm\quad
By noting that $\overline p_t=\lambda \overline {p(1-p)}(t) $, applying the variation-of-constants formula to the $p$-equation in \dref{1.1}
 yields
\begin{equation*}
\begin{array}{ll}\label{3.33}
p(\cdot,t)-\overline p(t)=&
e^{t\Delta}(p(\cdot,0)-\overline p(0))-\displaystyle\alpha\int^t_{0} e^{(t-s)\Delta} ( \displaystyle\frac p {1+c} c_x)_x\\
& -\rho\displaystyle\int^t_{0} e^{(t-s)\Delta} (  p w_x)_x
+\lambda \displaystyle\int^t_{0} e^{(t-s)\Delta} (p(1-p)-\overline{p(1-p)}).
\end{array}
\end{equation*}
Together with (3.4) and Lemma \ref{lemma315} and Lemma \ref{lemma316}, this gives
\begin{equation}
\begin{array}{rl}\label{3.36}
& \|p(\cdot,t)-\overline p(t)\|_{L^\infty(\Omega)}\\[2mm]
\leq & \|e^{t\Delta}(p(\cdot,0)-\overline p(0))\|_{L^\infty(\Omega)}+
\alpha\displaystyle\int^t_{0}\| e^{(t-s)\Delta} ( \displaystyle\frac p {1+c} c_x)_x\|_{L^\infty(\Omega)}\\[2mm]
&+\rho\displaystyle\int^t_{0} \|e^{(t-s)\Delta} (  p w_x)_x\|_{L^\infty(\Omega)}
+\lambda \displaystyle\int^t_{0} \|e^{(t-s)\Delta} (p(1-p)-\overline{p(1-p)})\|_{L^\infty(\Omega)}\\[2mm]
\leq& k_1e^{-\lambda_1 t}
 \|p(\cdot,0)-\overline p(0)\|_{L^\infty(\Omega)}+
 c_1\displaystyle\int^t_{0}(1+(t-s)^{-\frac 34}) e^{-\lambda_1(t-s)}
  \| w_x\|_{L^2(\Omega)}\\[2mm]
&+ c_1\displaystyle\int^t_{0}(1+(t-s)^{-\frac 34}) e^{-\lambda_1(t-s)} \| c_x\|_{L^2(\Omega)}\\[2mm]
&+ c_1\displaystyle\int^t_{0}(1+(t-s)^{-\frac 34}) e^{-\lambda_1(t-s)} \|p(1-p)-\overline{p(1-p)}\|_{L^2(\Omega)}\\[2mm]
\leq&k_1e^{-\lambda_1 t}
 \|p(\cdot,0)-\overline p(0)\|_{L^\infty(\Omega)}+
 c_2(\epsilon)\displaystyle\int^t_{0}(1+(t-s)^{-\frac 34}) e^{-\lambda_1(t-s)}e^{-\gamma(1-\epsilon)s}\\[2mm]
 &+c_2(\epsilon)\displaystyle\int^t_{0}(1+(t-s)^{-\frac 34}) e^{-\lambda_1(t-s)}e^{-(1-\epsilon)s}\\[2mm]
 &+ c_1\displaystyle\int^t_{0}(1+(t-s)^{-\frac 34}) e^{-\lambda_1(t-s)} \|p(1-p)-\overline{p(1-p)}\|_{L^2(\Omega)}.
\end{array}
\end{equation}
It is observed that
$$\|p(1-p)-\overline{p(1-p)}\|^2_{L^2(\Omega)}=\|p(1-p)\|^2_{L^2(\Omega)}-|\Omega||\overline{p(1-p)}|^2\leq \|p(1-p)\|^2_{L^2(\Omega)}.
$$
Hence from \eqref{3.33} and Lemma \ref{lemma315},  it follows that
\begin{equation}
\begin{array}{rl}\label{3.37}
\|p(\cdot,t)-\overline p(t)\|_{L^\infty(\Omega)}\leq & k_1e^{-\lambda_1 t}
 \|p(\cdot,0)-\overline p(0)\|_{L^\infty(\Omega)}+
 c_2(\epsilon)\displaystyle\int^t_{0}(1+(t-s)^{-\frac 34}) e^{-\lambda_1(t-s)}e^{-\gamma(1-\epsilon)s}\\
 &+c_2(\epsilon)\displaystyle\int^t_{0}(1+(t-s)^{-\frac 34}) e^{-\lambda_1(t-s)}e^{-(1-\epsilon)s}\\
 &+ c_3(\epsilon)\displaystyle\int^t_{0}(1+(t-s)^{-\frac 34}) e^{-\lambda_1(t-s)} e^{-\min\{1,\gamma,\lambda\}(1-\epsilon) s}\\
 \leq & c_4(\epsilon) e^{-\min\{\lambda_1,1,\gamma,\lambda\}(1-\epsilon) t},
\end{array}
\end{equation}
which implies \eqref{3.35}.

\vskip 3mm

\begin{lemma}\label{lemma318} Let $w_0>1-\displaystyle\frac 1 \rho$. Then for any  $\epsilon\in(0,\min\{\lambda_1,1,\gamma,\lambda\})$, there exists $C(\epsilon)>0$ such that
\begin{equation}\label{3.38}
\displaystyle |\overline p(t)-1| \leq C(\epsilon)e^{-(\min\{2\lambda_1,2,2\gamma,\lambda\}-\epsilon) t}.
\end{equation}
\end{lemma}
\it{Proof.}\rm\quad
We integrate the $p$-equation in the spatial variable over $\Omega$ to obtain
\begin{equation}\label{3.39}
\begin{array}{rl}
(\overline p-1)_t=&\lambda (\overline p-\overline p^2-\displaystyle\frac 1{|\Omega|}\int_\Omega(p-\overline p)^2)\\
=&-\lambda
\overline p(\overline p-1)-\displaystyle\frac \lambda{|\Omega|}\displaystyle \|p-\overline p\|^2_{L^2(\Omega)}.
\end{array}
\end{equation}
By \eqref{3.11}, 
 there exists $t_\epsilon>0$ such that
$
\overline p(t)\geq 1-\epsilon 
$
 for $t\geq t_\epsilon $. Hence by \eqref{3.35} and \eqref{3.39}, solving the differential equation entails
  \begin{equation*}\label{3.40}
\begin{array}{rl}
 |\overline p(t)-1|\leq & |\overline p(t_\varepsilon)-1|e^{-\lambda \int^t_{t_\varepsilon}\overline p(s)ds}+
 \displaystyle\frac \lambda{|\Omega|}\displaystyle \int^t_{t_\varepsilon}e^{-\lambda \int^t_{s}\overline p(\sigma)d\sigma}\|p(s)-\overline p(s)\|^2_{L^2(\Omega)}\\
\leq &|\overline p(t_\varepsilon)-1|e^{-\lambda(1-\epsilon)(t-t_\epsilon)} +
 \displaystyle\frac \lambda{|\Omega|}\displaystyle \int^t_{t_\varepsilon}e^{-\lambda (1-\epsilon)(t-s)}  \|p(s)-\overline p(s)\|^2_{L^2(\Omega)}\\[2mm]

\leq &|\overline p(t_\varepsilon)-1|e^{-\lambda(1-\epsilon)(t-t_\epsilon)} +
 c_1(\epsilon)\displaystyle \int^t_{0}e^{-\lambda (1-\epsilon)(t-s)}
e^{-2\min\{\lambda_1,1,\gamma,\lambda\}(1-\epsilon) s} \\[2mm]

\leq &|\overline p(t_\varepsilon)-1|e^{-\lambda(1-\epsilon)(t-t_\epsilon)} +
 c_2(\epsilon)\displaystyle
e^{-\min\{2\lambda_1,2,2\gamma,\lambda\}(1-\epsilon) t} \\
\leq &
c_3(\epsilon)\displaystyle
e^{-\min\{2\lambda_1,2,2\gamma,\lambda\}(1-\epsilon) t}
 \end{array}
\end{equation*}
for $t\geq t_\epsilon $, which proves \eqref{3.38}.

\vskip 3mm

\begin{lemma}\label{lemma318} Let $w_0>1-\displaystyle\frac 1 \rho$. Then for any $\epsilon\in(0,\min\{\lambda_1,1,\gamma,\lambda\})$, there exists $C(\epsilon)>0$ such that
\begin{equation}\label{3.41a}
\displaystyle \|p(\cdot,t)-1\|_{L^\infty(\Omega)} \leq C(\epsilon)e^{-(\min\{\lambda_1,1,\gamma,\lambda\}-\epsilon) t}.
\end{equation}
\end{lemma}
\it{Proof.}\rm\quad Combining the above two lemmas, we have
\begin{equation*}
\displaystyle \|p(\cdot,t)-1\|_{L^\infty(\Omega)} \leq
\displaystyle \|p(\cdot,t)-\overline p(t)\|_{L^\infty(\Omega)}+\displaystyle |\overline p(t)-1| \leq
C(\epsilon)e^{-(\min\{\lambda_1,1,\gamma,\lambda\}-\epsilon) t}.
\end{equation*}

\vskip 3mm

{\it Proof of  Theorem 1.2.}~\eqref{1.5} is the direct consequence of Lemma \ref{lemma38} in the previous subsection. As for \eqref{1.6}--\eqref{1.8}, we only need to collect (3.4), (3.26), (3.33) and (3.41).

\vskip3mm

{\bf Remark 3.2.}~ In comparison with \eqref{3.41a}, by \eqref{3.8} and \eqref{3.28}, we have $\sup_{t\geq 0}\int_\Omega | w_x(t)|^2\leq C$ and
$\sup_{t\geq 0}$ $\int_\Omega | z_x(t)|^2\leq C$ respectively, and thus $\sup_{t\geq 0}\|p(t)\|_{W^{1,2}(\Omega)}\leq C$. Hence an  interpolation by means of the Gagliardo--Nirenberg inequality in the one-dimensional setting  provides
 \begin{equation*}\label{3.41}
\displaystyle \|p(\cdot,t)-1\|_{L^\infty(\Omega)} \leq
\displaystyle \|p(\cdot,t)\|^{\frac 12}_{W^{1,2}(\Omega)}\displaystyle \|p(\cdot,t)-1\|^{\frac 12}_{L^2(\Omega)} \leq
C(\epsilon)e^{-(\frac 12\min\{1,\gamma,\lambda\}-\epsilon) t},
\end{equation*}
 where we have used (3.34).

 \section*{Acknowledgment}
  The authors are grateful to the referee for his illuminating comments.  This work is partially supported by the NUS AcRF grant
R-146-000-249-114 (PYHP) and by the NNSFC grant 11571363 (YW).

\end{document}